# Effective Prime Factors of the Given Even Number, Goldbach's Theorem and Goldbach's Conjecture

Song Linggen
(Dept. Materials Science, Fudan Univ., Shanghai (200433), P. R. China)

## Keywords:

Effective primes of the given even number greater than *6*;
Similarity theorem;
Two-part method;
Goldbach's theorem;
Goldbach's conjecture.

## Abstract:

Other than any odd prime whose factor is contained by the given even number $2N_m$ greater than 6, the odd primes $p_{(i, m)}$ within open interval $(1, 2N_m - 1)$ were defined as effective primes of $2N_m$.

Let    $2N_m - p_{(i, m)} = \alpha_{(i, m)}$    $2N_m - p_{(1, m)} p_{(i, m)} = \beta_{(i, m)}$

where,   Subscript *m* of $2N_m$ and $p_{(i, m)}$ represents that there are totally *m* effective primes of the given even number $2N_m$ and $p_{(i, m)}$ is the *i*-th minimum one of the *m* effective primes of $2N_m$, $1 < p_{(1, m)} < p_{(2, m)} < p_{(3, m)} < \ldots < p_{(m-1, m)} < p_{(m, m)} < 2N_m - 1$;
$p_{(1, m)}$ and $p_{(m, m)}$ is the minimum and maximum effective primes of $2N_m$, respectively;
$m \geqslant 2$ (to be verified).

Let symbol $D \langle\!\langle P /\!/ \alpha, \beta \rangle\!\rangle$ represent the distribution of prime factors *P* over integers $\alpha$ and $\beta$.

According to the similarity between $D \langle\!\langle p_{(j, M+1)} (j \neq M+1) /\!/ \alpha_{(j, M+1)} (j \neq M+1), \beta_{(j, M+1)} (j \neq M+1) \rangle\!\rangle$ when $m = M+1$ and $D \langle\!\langle p_{(i, M)} /\!/ \alpha_{(i, M)}, \beta_{(i, M)} \rangle\!\rangle$ when $m = M$ and the restrictions resulted from the relatively prime theorem on such distributions, by using two-part method, Goldbach's theorem was proved that whatever $m$ ($m \geqslant 2$) is, at least two of $\beta_{(i, m)}$ do not contain any factor of the effective primes $p_{(i, m)}$ of $2N_m$ and, consequently, at least two of $\alpha_{(i, m)}$ equals an effective prime of the given even number $2N_m$ greater than 6.

According to Goldbach's theorem, Goldbach's conjecture was proved.

## Introduction:

Other than any odd prime whose factor is contained by the given even number $2N_m$ greater than 6, the odd primes $p_{(i, m)}$ within open interval $(1, 2N_m - 1)$ were defined as effective primes of $2N_m$.

Let $\quad 2N_m - p_{(i, m)} = \alpha_{(i, m)} \qquad 2N_m - p_{(1, m)} p_{(i, m)} = \beta_{(i, m)}$

where, Subscript $m$ of $2N_m$ and $p_{(i, m)}$ represents that there are totally $m$ effective primes of the given even number $2N_m$ and $p_{(i, m)}$ is the $i$-th minimum one of the $m$ effective primes of $2N_m$, $1 < p_{(1, m)} < p_{(2, m)} < p_{(3, m)} < \ldots < p_{(m-1, m)} < p_{(m, m)} < 2N_m - 1$;
$p_{(1, m)}$ and $p_{(m, m)}$ is the minimum and maximum effective primes of $2N_m$, respectively;
$m \geqslant 2$ (to be verified).

Goldbach's conjecture states that any given even number $2N_m$ greater than 6 can be represented as a sum of two odd primes.

Goldbach's theorem to be proved states whatever $m$ ($m \geqslant 2$) is, at least two of $\beta_{(i, m)}$ do not contain any factor of the effective primes $p_{(i, m)}$ of $2N_m$ and, consequently, at least two of $\alpha_{(i, m)}$ equals an effective prime of the given even number $2N_m$ greater than 6, respectively.

Enlightened by Ben Green and Terence Tao[1], this paper attempted to verify Goldbach's theorem from which it can be derived that Goldbach's conjecture holds. .

## Verification of Goldbach's Theorem

*Effective products/integers of the given even number $2N_m$*
Products which contain and only contain the effective prime factor(s) of $2N_m$ were defined as effective products of $2N_m$. Effective primes and effective products of $2N_m$ were commonly named as effective integers of $2N_m$.

There are at least two different effective primes of any even number $2N_m$ greater than 6, because both integers $(N_m - k)$ and $(N_m + k)$ (where, $k = 1$ or 2 when $N$ is even or odd, respectively) are relatively prime to each other, and they contains at least one effective primes factor of $2N_m$, respectively, when $N_m$ is greater than 3. Therefore, $m \geq 2$.

*Restrictions on distributions of $m$ effective prime factors of $2N_m$ over $\alpha_{(i, m)}$ and $\beta_{(i, m)}$*
Because of relatively prime theorem, $D \langle\!\langle p_{(i, m)} // \alpha_{(i, m)}, \beta_{(i, m)} \rangle\!\rangle$, that is, the distributions of the $m$ effective prime factors $p_{(i, m)}$ of $2N_m$ over $\alpha_{(i, m)}$ and $\beta_{(i, m)}$ ($i = 1, 2, 3, \ldots, m$, respectively), are under the following restrictions:

*restr. 1*
Each of $a_{(i, m)}$ is an effective integer of $2N_m$ within open interval $(1, 2N_m - 1)$;

*restr. 2*
$a_{(a, m)}$ doesn't contain factor $p_{(a, m)}$, where, subscript *a* represents *1, 2, 3, …, m*, respectively.

*restr. 3*
$\beta_{(b, m)}$ doesn't contain factors $p_{(1, m)}$ and $p_{(b, m)}$, where, subscript *b* represents *1, 2, 3, …, m*, respectively.

*restr. 4*
There is no any common divisor between $a_{(c, m)}$ and $\beta_{(c, m)}$, where, subscript *c* represents *1, 2, 3, …, m*, respectively.

*restr. 5*
If $\beta_{(d, m)}$ doesn't contain any factor of the effective primes $p_{(i, m)}$ of $2N_m$, $1 \geqslant \beta_{(d, m)} > p_{(1, m)} p_{(d, m)}$ and $p_{(1, m)} p_{(m, m)} \geqslant p_{(1, m)} p_{(d, m)} \geqslant 2N_m - 1$

*restr. 6*
At most one of $a_{(i, m)}$ and $\beta_{(i, m)}$ may contain factor of the maximum effective prime $p_{(m, m)}$ of $2N_m$, respectively.

*restr. 7*
Because $p_{(1, m)} p_{(d, m)}$ is the minimum effective product containing factor $p_{(d, m)}$, there is no any effective product containing factor $p_{(d, m)}$ within open interval $(1, 2N_m-1)$ if $p_{(1, m)} p_{(d, m)} \geqslant 2N_m - 1$.

*Similarity theorem on the distributions of the effective prime factors*
Let $p_{(i,M)}$ represent the *M* effective primes of $2N_M$, $i = 1, 2, 3, …, M$, respectively, when $m = M$;
$p_{(j,M+1)}$ represent the *M+1* effective primes of $2N_{M+1}$, $j = 1, 2, 3, …, M, M+1$, respectively, when $m = M+1$;
symbol $D \langle\!\langle P // a, \beta \rangle\!\rangle$ represent the distribution of prime factors *P* over integers $a$ and $\beta$.

The similarity theorem states that:
Because of the similarity between $D \langle\!\langle p_{(j,M+1)}(j \neq M+1) // a_{(j,M+1)}(j \neq M+1), \beta_{(j,M+1)}(j \neq M+1) \rangle\!\rangle$ when $m = M+1$ and $D \langle\!\langle p_{(i, M)} // a_{(i,,M)}, \beta_{(i, M)} \rangle\!\rangle$ when $m = M$, if at least two of $\beta_{(i,,M)}$ do not contain any factor of the *M* effective primes $p_{(i, M)}$ of $2N_M$, at least two of $\beta_{(j, M+1)}(j \neq M+1)$ do not contain any factor of the *M* effective primes $p_{(j, M+1)}(j \neq M+1)$ of $2N_{M+1}$.

*The verification of Goldbach's theorem by using two-part method*
According to the restrictions and the similarity theorem mentioned, by using two-part method starting from $m = 2$, Goldbach's theorem was verified as shown below:

Part One
When $m = 2$.

There are totally two effective primes of the given even number $2N_2$ greater than 6, $p_{(1, 2)} < p_{(2, 2)}$, when $m = 2$ and

$$2N_2 - p_{(1, 2)} = \alpha_{(1, 2)} \qquad 2N_2 - p_{(1, 2)}p_{(1, 2)} = \beta_{(1, 2)}$$
$$2N_2 - p_{(2, 2)} = \alpha_{(2, 2)} \qquad 2N_2 - p_{(1, 2)}p_{(2, 2)} = \beta_{(2, 2)}$$

Observing $D \langle\!\langle p_{(i, 2)} // \alpha_{(i, 2)}, \beta_{(i, 2)} \rangle\!\rangle$ gave that restricted by <u>restr.s 1, 2, 3, 4,</u> both $\beta_{(1, 2)}$ and $\beta_{(2, 2)}$ do not contain any factor of effective primes $p_{(1, 2)}$ and $p_{(2, 2)}$ of $2N_2$ and, consequently, according to <u>restr.s 5, 7,</u> $\alpha_{(1, 2)}$ and $\alpha_{(2, 2)}$ equals $p_{(2, 2)}$ and $p_{(1, 2)}$, respectively.

Goldbach's theorem holds when $m = 2$.

Part Two
It was verified that if Goldbach's theorem holds when $m = M$, the theorem holds, too, when $m = M+1$ ($M \geq 2$) as shown below:

Comparing $D \langle\!\langle p_{(j, M+1)} // \alpha_{(j, M+1)} (j \neq M+1), \beta_{(j, M+1)}(j \neq M+1) \rangle\!\rangle$ with $D \langle\!\langle p_{(i, M)} // \alpha_{(i, M)}, \beta_{(i, M)} \rangle\!\rangle$ gave that one more factor, the maximum effective prime factor $p_{(M+1, M+1)}$ of $2N_{M+1}$, is involved in the first mentioned of two than the latter.

According to <u>restr. 6,</u> in terms of how the maximum effective prime factor $p_{(M+1, M+1)}$ of $2N_{M+1}$ involves in $\alpha_{(j, M+1)} (j \neq M+1)$ and $\beta_{(j, M+1)} (j \neq M+1)$, there are totally the following possible 4 cases in principle:

*Case 1*, None of $\alpha_{(j, M+1)} (j \neq M+1)$ and $\beta_{(j, M+1)} (j \neq M+1)$ contains factor $p_{(M+1, M+1)}$ of $2N_{M+1}$;

*Case 2*, Only $\alpha_{(a, M+1)}$ contains factor $p_{(M+1, M+1)}$ of $2N_{M+1}$;

*Case 3*, Only $\beta_{(b, M+1)}$ contains factor $p_{(M+1, M+1)}$ of $2N_{M+1}$;

*Case 4*, Both $\alpha_{(a, M+1)}$ and $\beta_{(b, M+1)}$ contain factor $p_{(M+1, M+1)}$ of $2N_{M+1}$.

where, subscripts *a* and *b* are two of *1, 2, 3, ..., M -1, M*.

It was verified that Goldbach's theorem holds in each logical case above when $m = M+1$ ($M \geq 2$) if Goldbach's theorem holds when $m = M$ as shown below:

<u>(1), When none of $\alpha_{(j, M+1)} (j \neq M+1)$ and $\beta_{(j, M+1)} (j \neq M+1)$ contains factor $p_{(M+1, M+1)}$:</u>
$D \langle\!\langle p_{(j, M+1)} // \alpha_{(j, M+1)} (j \neq M+1), \beta_{(j, M+1)}(j \neq M+1) \rangle\!\rangle$ is equivalent to $D \langle\!\langle p_{(j, M+1)} (j \neq M+1) //$

$a_{(j,M+1)}$ ($j \neq M+1$), $\beta_{(j, M+1)}$ ($j \neq M+1$)》 which is exactly similar to $D$ 《$p_{(i, M)} // a_{(i,,M)}, \beta_{(i, M)}$》 when none of $a_{(j,M+1)}$ ($j \neq M+1$) and $\beta_{(j,M+1)}$ ($j \neq M+1$) contains factor $p_{(M+1,M+1)}$.

Because of such similarity mentioned above, at least two of $\beta_{(j, M+1)}$ ($j \neq M+1$) do not contain any factor of the $M+1$ effective primes $p_{(j,M+1)}$ of $2N_{M+1}$ and, consequently, at least two of $a_{(j,M+1)}$ ($j \neq M+1$) equals an effective prime of the given even number $2N_{M+1}$, respectively, if Goldbach's theorem holds when $m = M$.

That is, at least two of $\beta_{(j, M+1)}$ do not contain any factor of the $M+1$ effective primes $p_{(j,M+1)}$ of $2N_{M+1}$ and, consequently, at least two of $a_{(j,M+1)}$ equals an effective prime of the given even number $2N_{M+1}$, respectively, if Goldbach's theorem holds when $m = M$.

Goldbach's theorem holds in *Case 1* when $m = M+1$ if the theorem holds when $m = M$.

*(2), When only $a_{(a, M+1)}$ contains factor $p_{(M+1, M+1)}$*

According to the similarity theorem, at least one of $\beta_{(j,M+1)}$ ($j \neq M+1$), $\beta_{(b,M+1)}$ ($b \neq a$), still doesn't contain any factor of the $M+1$ effective primes $p_{(j,M+1)}$ of $2N_{M+1}$ if Goldbach's theorem holds when $m = M$.

According to <u>restr.s 1, 5, 7,</u> $a_{(a, M+1)} = p_{(M+1, M+1)}$ when $\beta_{(b,M+1)}$ doesn't contain any factor of the $M+1$ effective primes $p_{(j,M+1)}$ of $2N_{M+1}$.

$a_{(M+1, M+1)} = p_{(a, M+1)}$ when $a_{(a, M+1)} = p_{(M+1, M+1)}$.

According to <u>restr. 5</u>, $\beta_{(M+1,M+1)}$ may doesn't contain any factor of the $M+1$ effective primes $p_{(j,M+1)}$ of $2N_{M+1}$, or may be a positive or negative effective integer of $2N_{M+1}$.

- If $\beta_{(M+1,M+1)}$ doesn't contain any factor of the $M+1$ effective primes $p_{(j,M+1)}$ of $2N_{M+1}$, at least two of $\beta_{(j,M+1)}$, $\beta_{(b,M+1)}$ and $\beta_{(M+1,M+1)}$, do not contain any factor of $M+1$ effective primes $p_{(j,M+1)}$ of $2N_{M+1}$ and at least $a_{(a, M+1)}$ and $a_{(M+1, M+1)}$ equals $p_{(M+1, M+1)}$ and $p_{(a, M+1)}$, respectively.

Goldbach's theorem holds in this sub-case if Goldbach's theorem holds when $m = M$.
.
- According to <u>restr. 4,</u> $\beta_{(M+1,M+1)}$ doesn't contain factor $p_{(a, M+1)}$ when $a_{(M+1, M+1)} = p_{(a, M+1)}$.

According to <u>restr. 3</u>, $\beta_{(M+1,M+1)}$ doesn't contain factor $p_{(M+1, M+1)}$.

There is no any common divisor between $\beta_{(M+1,M+1)}$ and $\beta_{(a,M+1)}$, because the sum, $\beta_{(a,M+1)} + \beta_{(M+1,M+1)} = (2 - p_{(1, M+1)})2N_{M+1}$ when $a_{(a, M+1)} = p_{(M+1, M+1)}$, doesn't contain any factor of the $M+1$ effective primes $p_{(j,M+1)}$ of $2N_{M+1}$.

If $\beta_{(M+1,M+1)}$ is a positive or negative effective integer of $2N_{M+1}$, because $\beta_{(M+1,M+1)}$ doesn't contain factors $p_{(a, M+1)}$ and $p_{(M+1, M+1)}$ and there is no any common divisor between $\beta_{(M+1,M+1)}$ and $\beta_{(a,M+1)}$ as verified above, $D \langle\!\langle p_{(j, M+1)}( j\neq M+1)// \alpha_{(j,M+1)}( \alpha_{(a, M+1)}$ was replaced by $\beta_{(M+1,M+1)}, j\neq M+1), \beta_{(j, M+1)} ( j\neq M+1)\rangle\!\rangle$ is similar to $D \langle\!\langle p_{(i, M)}// \alpha_{(i,,M)}, \beta_{(i, M)}\rangle\!\rangle$.

Because of such similarity mentioned above, at least two of $\beta_{(j, M+1)}( j\neq M+1)$ do not contain any factor of the effective primes $p_{(j,M+1)}$ of $2N_{M+1}$.

That is, at least two of $\beta_{(j, M+1)}$ do not contain any factor of the effective primes $p_{(j,M+1)}$ of $2N_{M+1}$, and at least $\alpha_{(a, M+1)}$ and $\alpha_{(M+1, M+1)}$ equals $p_{(M+1, M+1)}$ and $p_{(a, M+1)}$, respectively.

Goldbach's theorem holds in this sub-case if Goldbach's theorem holds when $m = M$.

Sub-summarizing: Goldbach's theorem holds in *Case 2* when $m = M+1$ if the Theorem holds when $m = M$.

*(3), When only $\beta_{(b, M+1)}$ contains factor $p_{(M+1, M+1)}$*

According to the similarity theorem, at least one of $\beta_{(j,M+1)} ( j\neq M+1), \beta_{(b',M+1)} ( b'\neq b)$ still doesn't contain any factor of the $M+1$ effective primes $p_{(j,M+1)}$ of $2N_{M+1}$ if Goldbach's theorem holds when $m = M$.

According to <u>restr.s 5, 7</u>, $\beta_{(b, M+1)}$ has to equal $p_{(M+1, M+1)}$ or $- p_{(M+1, M+1)}$ when at least $\beta_{(b',M+1)}$ doesn't contain any factor of the $M+1$ effective primes $p_{(j,M+1)}$ of $2N_{M+1}$.

According to <u>restr.s 3,</u> $\beta_{(M+1,M+1)}$ doesn't contain factor $p_{(M+1, M+1)}$.

- When $\beta_{(b,M+1)} = p_{(M+1, M+1)}$, $\alpha_{(M+1,M+1)} = p_{(1, M+1)}p_{(b, M+1)}$ and consequently, according to <u>restr.4</u>, $\beta_{(M+1,M+1)}$ doesn't contain factor $p_{(b, M+1)}$ when $\alpha_{(M+1, M+1)} = p_{(1, M+1)}p_{(b, M+1)}$.

  There is no any common divisor between $\beta_{(M+1,M+1)}$ and $\alpha_{(b,M+1)}$, because the sum, $\beta_{(M+1,M+1)} \alpha_{(b,M+1)}p_{(1, M+1)}{}^2 = p_{(1, M+1)}( p_{(1, M+1)} -1)2N_{M+1)}$ when $\beta_{(b, M+1)} = p_{(M+1, M+1)}$, doesn't contain any factor of the $M+1$ effective primes $p_{(j,M+1)}$ of $2N_{M+1}$.

  Because $\beta_{(M+1,M+1)}$ doesn't contain factors $p_{(b, M+1)}$ and $p_{(M+1, M+1)}$ and there is no any common divisor between $\beta_{(M+1,M+1)}$ and $\alpha_{(b,M+1)}$, $D \langle\!\langle p_{(j, M+1)}( j\neq M+1)// \alpha_{(j,M+1)}( j\neq M+1), \beta_{(j, M+1)} (\beta_{(b,M+1)}$ was replaced by $\beta_{(M+1,M+1)}, j\neq M+1)\rangle\!\rangle$ is similar to $D \langle\!\langle p_{(i, M)}// \alpha_{(i,,M)}, \beta_{(i, M)}\rangle\!\rangle$ when $\beta_{(b,M+1)} = p_{(M+1, M+1)}$.

  Because of such similarity mentioned above, at least two of $\beta_{(j, M+1)}( j\neq b)$ do not contain any factor of the effective primes $p_{(j,M+1)}$ of $2N_{M+1}$ and , consequently, at least two of $\alpha_{(j, M+1)}$ ( $j\neq M+1$) equals an effective prime of the given even number $2N_{M+1}$, respectively.

  That is, at least two of $\beta_{(j, M+1)}$ do not contain any factor of the effective primes $p_{(j,M+1)}$ of $2N_{M+1}$ and , consequently, at least two of $\alpha_{(j, M+1)}$ equals an effective prime of the given even number $2N_{M+1}$, respectively.

Goldbach's theorem holds in this sub-case if Goldbach's theorem holds when $m = M$.

- When $\beta_{(b, M+1)} = -p_{(M+1, M+1)}$, that is $2N_{M+1} + p_{(M+1, M+1)} = p_{(1, M+1)}p_{(b, M+1)}$ and, consequently, according to relatively prime theorem, $\alpha_{(M+1, M+1)}$ ( $= 2N_{M+1} - p_{(M+1, M+1)}$ ) doesn't contain factor $p_{(b, M+1)}$ when $\beta_{(b, M+1)} = -p_{(M+1, M+1)}$.

    There is no any common divisor between $\alpha_{(M+1, M+1)}$ and $\alpha_{(b, M+1)}$, because the difference, $\alpha_{(b, M+1)}p_{(1, M+1)} - \alpha_{(M+1, M+1)} = (p_{(1, M+1)} - 2)2N_{M+1}$ when $\beta_{(b, M+1)} = -p_{(M+1, M+1)}$, doesn't contain any factor of the $M+1$ effective primes $p_{(j, M+1)}$ of $2N_{M+1}$.

    Because $\alpha_{(M+1, M+1)}$ doesn't contain factors $p_{(b, M+1)}$ and $p_{(M+1, M+1)}$ and there is no any common divisor between $\alpha_{(M+1, M+1)}$ and $\alpha_{(b, M+1)}$, $D$ 《$p_{(j, M+1)}( j \neq M+1)$// $\alpha_{(j, M+1)}( j \neq M+1)$, $\beta_{(j, M+1)}$ ($\beta_{(b, M+1)}$ was replaced by $\alpha_{(M+1, M+1)}$, $j \neq M+1$)》 is similar to $D$ 《$p_{(i, M)}$// $\alpha_{(i, M)}$, $\beta_{(i, M)}$》 when $\beta_{(b, M+1)} = -p_{(M+1, M+1)}$.

    Because of such similarity mentioned above, at least two of $\beta_{(j, M+1)}( j \neq M+1, j \neq b)$ do not contain any factor of the effective primes $p_{(j, M+1)}$ of $2N_{M+1}$ and , consequently, at least two of $\alpha_{(j, M+1)}( j \neq M+1)$ equals an effective prime of the given even number $2N_{M+1}$, respectively.

    That is, at least two of $\beta_{(j, M+1)}$ do not contain any factor of the effective primes $p_{(j, M+1)}$ of $2N_{M+1}$ and , consequently, at least two of $\alpha_{(j, M+1)}$ equals an effective prime of the given even number $2N_{M+1}$, respectively.

    Goldbach's theorem holds in this sub-case if Goldbach's theorem holds when $m = M$.

Sub-summarizing: Goldbach's theorem holds in *Case 3* when $m = M+1$ if the Theorem holds when $m = M$.

*(4),* When both $\alpha_{(a, M+1)}$ and $\beta_{(b, M+1)}$ contain factor $p_{(M+1, M+1)}$ of $2N_{M+1}$.
According to the similarity theorem and <u>restr. 5</u>, the absolute value of $\alpha_{(a, M+1)}$ and $\beta_{(b, M+1)}$ has to be an exact power of factor $p_{(M+1, M+1)}$, respectively, when both $\alpha_{(a, M+1)}$ and $\beta_{(b, M+1)}$ contain factor $p_{(M+1, M+1)}$ of $2N_{M+1}$.

Because the absolute value of the difference between $\alpha_{(a, M+1)}$ and $\beta_{(b, M+1)}$ has to be less than $p_{(M+1, M+1)}(p_{(M+1, M+1)} - 1)$, according to <u>restr. 1</u>, $\alpha_{(a, M+1)} = p_{(M+1, M+1)}$ and $\beta_{(b, M+1)} = -p_{(M+1, M+1)}$ when both $\alpha_{(a, M+1)}$ and $\beta_{(b, M+1)}$ contain factor $p_{(M+1, M+1)}$ of $2N_{M+1}$.

By the same logics shown in *(2)* and *(3)* above, $D$ 《$p_{(j, M+1)} ( j \neq M+1)$// $\alpha_{(j, M+1)}( \alpha_{(a, M+1)}$ was replaced by $\beta_{(M+1, M+1)}$ ), $j \neq M+1$), $\beta_{(j, M+1)}( \beta_{(b, M+1)}$ was replaced by $\alpha_{(M+1, M+1)}$), $j \neq M+1$》 is similar to $D$ 《$p_{(i, M)}$// $\alpha_{(i, M)}$, $\beta_{(i, M)}$》 .

Because of such similarity mentioned above, by the same logic shown in *(2)* and *(3)* above Goldbach's theorem holds in *Case 4* when $m = M+1$ if the Theorem holds when $m = M$.

Summarizing above, if Goldbach's Theorem holds when $m = M$, the Theorem holds, too, when $m = M+1$.

Because Goldbach's theorem holds when $m = 2$ as verified in Part One, according to the result proved in Part Two, whatever $m$ $(m \geq 2)$ is, at least two of $\beta_{(i, m)}$ do not contain any factor of the effective primes $p_{(i, m)}$ of $2N_m$ and, consequently, at least two of $\alpha_{(i, m)}$ equals an effective prime of the given even number $2N_m$ greater than 6, respectively.

Goldbach's theorem holds.

## Verification of Goldbach's conjecture

Any even number greater than 6 is corresponding to a definite $m$ $(m \geq 2)$.

According to Goldbach's theorem, any given even number $2N_m$ greater than 6 can be represented as a sum of two odd primes.

## References

[1] Ben Green and Terence Tao, arXiv:math/0606088

**The end**